\documentclass[12pt]{amsart}
\usepackage{amssymb,amsmath,amsfonts}
\usepackage[all]{xy}

\begin{document}
\bibliographystyle{plain}

\def\secteqn{
\let\sectio\section%
\renewcommand{\section}{\sectioneqn\sectio }%
\newcommand{\sectioneqn}{\setcounter{equation}{0}%
  \renewcommand{\theequation}{\arabic{section}.\arabic{equation}}}}

\def\subsecteqn{
\let\subsectio\subsection%
\def\subsection{\subsectioneqn \subsectio}%
\def\subsectioneqn{\setcounter{equation}{0}%
  \def\theequation{\arabic{section}.\arabic{subsection}.\arabic{equation}}}
}

\def\RBimod{\mbox{$R$-\sffamily Bimod}}
\def\Rmod{\mbox{$R$-\sffamily mod}}
\def\modR{\mbox{\sffamily mod-$R$}}
\def\ModR{\mbox{\sffamily Mod-$R$}}
\def\Lcomod{\mbox{$L$-\sffamily comod}}
\def\comodL{\mbox{\sffamily comod-$L$}}
\def\comod{\mbox{\sffamily comod}}
\def\otr{\ot_R}
\def\RRbimod{\mbox{$R-R$-\sffamily bimod}}
\def\RRBimod{\mbox{$R-R$-\sffamily Bimod}}
\def\nablaop{\nabla^{\rm op}}

\def\sqot{\boxtimes}
\def\tsot{\hspace{.5ex}{}_\tau\hspace{-.5ex}\ot^\si}
\def\ttot{\hspace{.5ex}{}_\tau\hspace{-.5ex}\ot^\tau}
\def\stot{\hspace{.5ex}{}_\si\hspace{-.5ex}\ot^\tau}
\def\stdot{\hspace{.5ex}{}^\si\hspace{-.5ex}\ot_\tau}
\def\ssot{\hspace{.5ex}{}^\si\hspace{-1ex}\ot_\si}
\def\ttdot{\hspace{.5ex}{}^\tau\hspace{-1ex}\ot_\tau}
\def\ssuot{\hspace{.5ex}{}_\si\hspace{-.5ex}\ot^\si}

\def\rd#1{{}^*\hspace{-.5ex}#1}

\def\cat#1{{\mathfrak #1}}
\def\fun#1{{\mathcal #1}}

\def\labelenumi{(\roman{enumi})}

\def\example{\  \newline  \noindent {\sc Example. }}
\def\remark{\ \newline  \noindent {\sc Remark. }}
\def\definition{\ \newline  \noindent {\sc Definition. }}

\def\enlarge{\rule[2.1ex]{0cm}{0cm}}
\def\Enlarge{\rule[2.3ex]{0cm}{0cm}}
\def\largebr#1{\left(\enlarge #1\right)}
\def\Largebr#1{\left(\Enlarge #1\right)}

\def\lora{\longrightarrow}
\def\ot{\otimes}
\def\opl{\oplus}
\def\loma{\longmapsto}
\def\si{\sigma}

\theoremstyle{plain}
\newtheorem{thm}{Theorem}[subsection]
\newtheorem{lem}[thm]{Lemma}
\newtheorem{cor}[thm]{Corollary}
\newtheorem{pro}[thm]{Proposition}

\def\Ker{\mbox{\rm Ker\hspace{.1ex}}}
\def\Im{\mbox{\rm Im\hspace{.1ex}}}
\def\coend{\mbox{\sffamily Coend\hspace{.1ex}}}
\def\Mor{\mbox{\rm Mor\hspace{.1ex}}}
\def\Hom{\mbox{\rm Hom\hspace{.1ex}}}
\def\Vec{\mbox{\rm Vec }}
\def\End{\mbox{\rm End\hspace{.1ex}}}

\def\R{\overline{R}}

\def\eee{\rule{.75ex}{1.5ex}\vskip.5ex}

\def\proof{{\it Proof.\ }}

\newcommand{\va}{\varepsilon}
\newcommand{\var}{\varepsilon}
\newcommand{\Mod}{\mbox{Mod}}
\newcommand{\Comod}{\mbox{Comod}}
\def\lam{{\lambda}}

\renewcommand{\hom}{\underline{\mbox{\rm hom}}}
\newcommand{\Nat}{\mbox{\rm Nat}}

\def\rref#1{(\ref{#1})}

\def\db{{\mathchoice{\mbox{\rm db}}
                    {\mbox{\rm db}}
                    {\mbox{\scriptsize\rm db}}
                    {\mbox{\tiny\rm db}} }}

\def\ev{{\mathchoice{\mbox{\rm ev}}
                    {\mbox{\rm ev}}
                    {\mbox{\scriptsize\rm ev}}
                    {\mbox{\tiny\rm ev}} }}

\def\id{{\mathchoice{\mbox{\rm id}}
                    {\mbox{\rm id}}
                    {\mbox{\scriptsize\rm id}}
                    {\mbox{\tiny\rm id}} }}

\def\op{{\mathchoice{\mbox{\rm op}}
                    {\mbox{\rm op}}
                    {\mbox{\scriptsize\rm op}}
                    {\mbox{\tiny\rm op}} }}

\author{Ph\`ung H\^o  Hai}
\title{Tannaka-Krein duality for Hopf algebroids}
\dedicatory{Institute of Mathematics \\ P.O.Box 631, 10000 Bo Ho, Hanoi, Vietnam
\\
phung@math.ac.vn\\[2ex]
{$\mathfrak {Herrn\ Prof.\ B.\hspace{.1ex} Pareigis\ zum\ 65\ Geburtstag\ gewidmet}$}}
 \thanks{Current address: Department of Mathematics, University of Duisburg-Essen, 45117 Essen,
Germany}

\begin{abstract} We show that a Hopf algebroid can be reconstructed from
 a monoidal functor from a monoidal category into the category of rigid
bimodules over a ring. We study the equivalence between the original category
and the category of comodules over the reconstructed Hopf algebroid.
\end{abstract}
\maketitle
\secteqn

\section*{Introduction}
The Tannaka-Krein duality asserts that a compact group can be uniquely determined by the category of its
finite dimensional unitary representations. Many efforts have been made in generalizing this
result, which also advance the development of many branches of mathematics, such as $C^*$-algebras,
harmonic analysis, algebraic geometry. Especially, Tannaka-Krein duality was also one
of the sources of quantum groups.

The algebraic version of this theory was suggested by
A. Grothendieck and developed by Saavedra, Deligne
\cite{saavedra,deligne90}.
An important result of the algebraic Tannaka-Krein theory is
a theorem of Deligne, developing Saavedra's ideas. It states
that there is a dictionary between tensor categories over a field
$k$ together with an exact tensor functor (fiber functor) to the
category of quasi-coherent sheaves over a $k$-scheme $S$ and transitive
groupoids over $S$.

The proof of Tannaka-Krein duality is divided into two parts, the
(re)construction theorem which aims to reconstruct the group
from the category of its representations and the representation theorem which
aims to prove the equivalence between the original category
and the category of representations of the reconstructed group.

The idea of the reconstruction theorem was one of the motivations for
Quantum Groups. From this point of view
a rigid monoidal category (without a symmetry)
corresponds to a quantum group. In fact, one can
construct from a rigid monoidal category together with a monoidal
functor into the category of finite-dimensional vector spaces (over $k$)
a $k$-Hopf algebra, which is understood as the ``function algebra"
over a quantum group. This idea was first proposed in the  work  \cite{lyu1} of
Lyubashenko. Tannaka-Krein duality for compact quantum
groups was proved by Woronowicz \cite{woron2}. Majid obtained
the reconstruction theorem in a more general setting of a monoidal
category and a monoidal functor in to another braided monoidal
category, see \cite{majid_book} and references therein.

In this more general setting, a Hopf algebra (in a braided monoidal category)
cannot be reconstructed from its representation (comodule) category. Usually,
the reconstructed Hopf algebra is bigger, see for instance \cite{par1}.
 Lyubashenko \cite{lyu2} suggests reconstructing the Hopf algebra not
lying in the target category but rather in its tensor square.
 McCrudden \cite{crudden1} also generalizes the duality 
 to the setting of higher categories.
  
An important ingredient of
 the Tannaka-Krein duality is the fiber functor. One might ask,
for what kind of monoidal categories there exist such functors.
An answer to this question can be called an embedding theorem.
For a tensor category over a field $k$ of characteristic 0, P. Deligne \cite{deligne90}
  gave an interesting criterion in terms of the categorical 
  dimension. A parallel  result for $C^*$ tensor
categories was given by Doplicher and Roberts \cite{dr89}.

In our previous work \cite{ph00} an embedding theorem for arbitrary (rigid)
monoidal categories was given, but the embedding goes into the category
of bimodules over a ring. This result raises the problem
of Tannaka-Krein duality for functors with target the category of
bimodules over a ring. However, according to Schauenburg \cite{schauen1},
 it is generally impossible to construct a braiding in a bimodule category, so
one cannot apply Majid results to reconstruct a Hopf algebra in this category.

It turns out that one can reconstruct from the above data 
a Hopf algebroid in the sense of  Takeuchi,
 Lu,  and Schauenburg \cite{takeuchi1,takeuchi2,lu1,schauen2}. This construction
 is in a sense analogous to those of Lyubashenko \cite{lyu2} and Mccrudden \cite{crudden1}.

Combining the Tannaka-Krein duality done here and the embedding theorem
of \cite{ph00}, we can realize a rigid category as the comodule category over a Hopf
algebroid defined over a certain ring  (Corollary \ref{cor_repr}).

The paper is constructed as follows. In Section 1, we  recall the notion
of Hopf algebroids defined over a ring. We define bialgebroid as a monoidal
 object in the monoidal category of coalgebroids. Next, we recall a notion
of antipode and Hopf algebroids and prove some basic facts on dual comodules
over Hopf algebroids. This is the more technically difficult part of the work. 
In fact there are at least two
definitions of antipode on a bialgebroid \cite{lu1,schauen2}. In \cite{lu1} the definition
of the antipode  more or less imitates the usual antipode, in \cite{schauen2} the
antipode is defined as a condition for the coincidence 
of internal hom-functors in the category
of modules and in the underlying category of $R$-bimodules. 
Our motivation for the antipode
is the condition for the existence of dual comodules over a bialgebroid. 
It is somewhat unexpected that the antipode introduced by Schauenburg in \cite{schauen2} while
studying duals of modules over a bialgebroid fits well into our frame-work.

In Section 2 we prove the Tannaka-Krein duality for Hopf algebroids.
Some  embedding and reconstruction results were also 
obtained by Hayashi \cite{haya960,haya96}  for face algebras, which
were shown by Schauenburg to be a special case of Hopf algebroids.
Our result here is a generalization of Hayashi's result.
\section{Hopf algebroids and its comodules}
Except for some results in Subsections \ref{sec1.3} and
\ref{sec1.7}, the materials of this section are
known, they can be found in \cite{takeuchi2,lu1,xu1,schauen2}.

We first review some basic notions of rings, corings
over an associative ring. In Section \ref{sec1.2} we recall the notion of 
coalgebroids which
was first studied by Takeuchi \cite{takeuchi2}.
In Subsections \ref{sec1.3} and \ref{sec1.4}
we study comodules over a coalgebroid and prove some lemmas
which will be needed in the sequel. In Subsection \ref{sec1.5} we recall the notion
of bialgebroids in the sense of \cite{takeuchi2,lu1}. In \ref{sec1.6} we define the
tensor product of two comodules over a bialgebroid and in \ref{sec1.7} the dual
to a comodule.

\subsection{$R$-rings and $R$-corings}\label{sec1.1}
Let us fix a commutative ring $k$. Throughout this paper, we will be working
in the category of $k$-modules,
in other words, we shall assume that everything is $k$-linear.

Let $R$ be an algebra over $k$, which will usually be fixed. Object of our study is the
category  $\RBimod$ of $R$-bimodules. In this category, there is a monoidal structure
with the tensor product being the  usual tensor product over $R$. This tensor product
is closed in the sense that there exist the right adjoint functors
to the functors $M\otr -$ and $-\otr M$ for all $R$-bimodules $M$, given by $\Hom_R(M,-)$
and ${}_R\Hom(M,-)$, where $\Hom_R(-,-)$ (resp. ${}_R\Hom(-,-))$
denotes the set of $R$-linear maps
with respect to the right (resp. left) actions of $R$.

Having the monoidal structure on $\RBimod$ we define $R$-rings
 and $R$-corings as monoids and comonoids in this category.
The data for an $R$-ring consist of
 an $R-R$-linear map $m:A\otr A\lora A$, called product, and an
 $R-R$-linear map $u:R \lora A$ called unit, satisfying
  the usual associativity
 and unity properties.
Set $1_A:=u(1_R)$ and denote $m(a \ot b)$ by $a \cdot b$,
then $A$ is a $k$-algebra in the usual sense.
We notice that if $R$ is commutative and $A$ is an algebra over 
$R$ in the usual sense
then it is an $R$-ring in our sense but the converse is not true 
since the image of
$R$ under $u$ is generally not in the center of $A$.
In fact, any (associative) $k$-algebra homomorphism
$R\lora A$ induces a structure of $R$-ring over $A$.

$R$-corings are defined in the dual way. A structure of $R$-coring
 over an $R$-bimodule $C$ consists of an $R-R$-linear map
$\Delta:C\lora C\otr C$ called coproduct and an $R-R$-linear
map $\varepsilon:C\lora R$, called counit, satisfying the usual
coassociativity and counity axioms.
We shall use Sweedler's notation for denoting the coproduct:
$$\Delta(a)=\sum_{(a)} a_{(1)}\ot a_{(2)}.$$
 A right $C$-comodule is a 
right $R$-module $M$ equipped with an $R$-linear
 coaction $\delta:M\lora M\otr C,\delta(m)=\sum_{(m)}m_{(0)}\ot m_{(1)},$
  satisfying
$$ \sum_{(m)}\delta(m_{(0)})\ot m_{(1)}=
\sum_{(m)}m_{(0)}\ot \Delta(m_{(1)})\quad\mbox{and}
\quad \sum_{(m)}m_{(0)}\ot \va(m_{(1)})=m$$
(having in mind the identification $M\otr R\cong M$). 
Notice that  the $R$-linearity of $\delta$ means
$\delta(mr)=\sum_{(m)}m_{(0)}\ot m_{(1)}r.$

\subsection{$R$-Coalgebroids}\label{sec1.2}
We consider in this subsection the category $\RRBimod$ of double $R$-bimodules.
That is, we will have two $R$-bimodule structures on a $k$-module, which commute with
each other. To distinguish the two structures we will denote the first one by $\si$
and the second one by $\tau$. Thus, we have four actions:
\begin{eqnarray*}\begin{array}{l@;l@,l@;l@,}
R\ot_k M\lora M&\quad  r\ot M\loma \si(r)m&\quad M\ot_k R\lora M&\quad  m\ot r\loma m\si(r)\\
R\ot_k M\lora M& \quad r\ot M\loma \tau(r)m&\quad M\ot_k R\lora M&\quad
m\ot r\loma m\tau(r)\end{array}\end{eqnarray*}

There are several possibilities to take tensor products over $R$ of $R-R$-bimodules.
 We use the notation $M\tsot N$ for the tensor product with respect to the right action
of $M$ by $\tau$ and the left action on $N$ by $\si$, i.e.,
$$ M\tsot N:=M\ot_kN\left/\Largebr{m\tau(r)\ot n=m\ot \si(r)n}\right. .$$
Here, the letters $\si$ and $\tau$  on the two sides of the tensor sign denote
 correspondingly the actions taken in the
definition of the tensor product.
Other tensor products will be denoted in a similar way. The rule for notation is
that the left action will be placed in the upper place
 and the right action will be placed in the lower place on the two sides of the tensor sign.

For the tensor product $M\tsot N$ we specify the following actions to make it
an object in $\RRBimod$:
$$\begin{array}{ll}
\tau(a)(h\ot k)=h\ot \tau(a)k,& (h\ot k)\tau(a)=h\ot k\tau(a),\\
\si(a)(h\ot k)=\si(a)h\ot k,& (h\ot k)\si(a)=h\si(a)\ot k.\end{array}$$
Here we adopt the convention that the action of $R$ has preference over the tensor product.

 \definition\cite{takeuchi2} An $R$-{\em  coalgebroid} is an $R-R$-bimodule $L$ equipped with
 $k$-linear maps $\Delta:L\lora L\tsot L$, called coproduct, and $\va:L\lora R$,
called counit, satisfying the following conditions:
\begin{enumerate}
\item $\Delta$ is a morphism in $\RRbimod$ and (the coassociativity):
$$ (\id_L\tsot \Delta)\Delta=(\Delta \tsot \id_L)\Delta,$$
\item $\va$ satisfies (the linearity with respect to the actions of $R$)
$$  \va(\si(a)h\tau(b))=a\va(h)b,$$
and (the counity)
 $$ (\va\tsot \id_L)\Delta=(\id_L\tsot \va)\Delta=\id_L.$$
\item[(iii)]Moreover, $\va$ satisfies the following condition
$$\va(\tau(a)h)=\va(h\si(a)).$$
\end{enumerate}
Note that, by definition, $\va$ is not necessarily a morphism of $R-R$-bimodules.\bigskip

We shall use Sweedler's notation for the coproduct: $\Delta(h)=\sum_{(h)}h_{(1)}\ot h_{(2)}.$
 The linearity of $\Delta$ now reads:
\begin{equation}\label{eq9} \Delta(\tau(a)\si(b)h\tau(c)\si(d))=\sum_{(h)}\si(b)h_{(1)}\si(d)\ot \tau(a)h_{(2)}\tau(c).
\end{equation}
Analogously, the counity condition has the following form
\begin{equation}\label{eq10} \sum_{(h)}\si\va(h_{(1)})h_{(2)}=\sum_{(h)}h_{(1)}\tau\va(h_{(2)})=h.
\end{equation}
Combining these equations, we have the following identities
\begin{equation}\label{eq11} h\si(a)=\sum_{(h)}\si\va\left(h_{(1)}\si(a)\right)h_{(2)};\quad
\tau(a)h=\sum_{(h)}h_{(1)}\tau\va\left(\tau(a)h_{(2)})\right),\end{equation}
whence
\begin{equation}\label{eq12} \va\left(\rule[2.1ex]{0cm}{0cm}\tau(a)h\si(b)\right)=
\sum_{(h)}\va(h_{(1)}\si(b))\va\left(\tau(a)h_{(2)}\right).\end{equation}

\remark The condition in (iii) can be replaced by the following (cf. \cite[\S3]{takeuchi2})
\begin{equation}\label{eq13}\sum_{(h)}\tau(a)h_{(1)}\tsot h_{(2)}=\sum_{(h)}h_{(1)}\tsot h_{(2)}\si(a).\end{equation}
This condition is equivalent to the existence
of an anchor as in \cite{xu1,lu1}. In fact, the algebra
$\End_k(R)$ has a structure of an $R-R$-bimodule, we specify it as follows:
$$ \largebr{\si(a)f\si(b)}(c):=af(bc),\quad \largebr{\tau(a)f\tau(b)}(c)=f(ca)b.$$
Then, the counity induces a morphism of $R-R$-bimodules $\eta:L\lora \End_k(R)$,
given by $\eta(h)(a):=\va(\tau(a)h)=\va(h\si(a)).$ The map 
$\eta$ is called an anchor \cite{xu1,lu1} (this map is generally different from a map,
also denoted by $\eta$, introduced in \cite[\S3]{takeuchi2}).

\subsection{An example} \label{sec_example}
Let $M$ be a right $R$-module. 
Then $M^*$ is a left $R$-module with the action given by
$$(r\varphi)(m):=r(\varphi(m));\quad r\in R, m\in M,\varphi\in M^*.$$
If $M$ is finitely generated (f.g.) projective, $M$ is a direct summand of $R^{\oplus d}$ considered as right module
over $R$, then $M^*$ is a direct summand of $R^{\oplus d}$ considered as left module over $R$.
 Further, if we fix a generating set $m_1,m_2,\ldots,m_d$ given by the projection from $R^{\oplus d}$ then we
can find a generating set $\varphi^1,\varphi^2,\ldots,\varphi^d$ for $M^*$ such that for any 
$m\in M$, the following equation holds true
\begin{equation}\label{eq1.6}m=\sum^d_{i=1}m_i\varphi^i(m).
\end{equation}
We define the following map
\begin{eqnarray}\label{eq1.7}&& \ev_{k,M}:M^*\otimes_kM\to R; \quad\varphi\otimes m
\mapsto \varphi(m),\\
\label{eq1.8} &&\db_{k,M}:k\to M\otimes_RM^*; 1\mapsto \sum_im_i\otimes \varphi^i.
\end{eqnarray}
Notice that $\ev_{k,M}$ is a morphism of $R$-bimodules:
$\ev_{k,M}(r\varphi\otimes ms)=(r\varphi)(ms)=r\varphi(m)s$. The equation in 
\rref{eq1.6} implies the following relations for $\ev=\ev_{k,M}$ and $\db=\db_{k,M}$
\begin{equation}\label{eq1.9}
(\ev\otimes_R\id_{M^*})(\id_{M^*}\otimes_k\db)=\id_{M^*};
\quad (\id_M\otimes_R\ev)(\db\otimes_k\id_{M})=\id_M.\end{equation}

Conversely, if there exists to a right $R$-module $M$ a left $R$-module $M^\vee$ and
morphisms $\ev:M^\vee\otimes_k M\to R$ and $\db:k\to M\otimes_RM^\vee$, satisfying the
identities in \rref{eq1.9} then $R$ is f.g. projective. Indeed, we have 
by means of \rref{eq1.9} the following natural isomorphism
$$ \Hom_R(P\otimes_k M,N)\cong \Hom_k(P,N\otimes_RM^\vee);\quad
f\mapsto (f\otimes_R\id_{M^\vee})(\id_P\otimes_k\db_M).$$
From the canonical isomorphism $\Hom_R(P\otimes_kM,N)\cong\Hom_k(P,
\Hom_R(M,N))$, we deduce a functorial isomorphism
$$N\otimes_RM^\vee\cong\Hom_R(M,N).$$
Since the functor $-\otimes_RM^\vee$ is right exact, $M$ is projective.
Setting $N=R$ in the isomorphism above 
we obtain isomorphism $M^\vee\cong\Hom_R(M,R)=M^*$, by means of which
the map $\ev$ is given by $\ev(\varphi\otimes m)=\varphi(m)$.
For $N=M$, the identity map $\id_M$ corresponds to the element
$\db(1)=\sum_{i=1}^dm_i\otimes \varphi^i$, with the property
$m=\sum^d_{i=1}m_i\varphi^i(m)$, for all $m\in M$. Hence $\{m_i\}$ generate $M$ and
$\{\varphi^i\}$ generate $M^*$.
We call the pair $\{m_i\}$, $\{\varphi^i\}$ \emph{dual bases} with respect to
$\db=\db_{k,M}$. In particular we have proved:
\begin{lem}
 Let $M$ be an f.g. projective right $R$-module. 
 Denote the action of $R$ on $M$
by $\tau$ and the one on $M^*$ by $\sigma$. 
Then $M^*\otimes_kM$ is an $R$-bimodule and $\ev_{k,M}$ is an $R$-bimodule homomorphism.
\end{lem}

Assume now that $M$ is an $R$-bimodule which is f.g. projective
as an $R$-module. The the left action $R$ on $M$ induces a right action of $R$ on
$M^*=\Hom_R(M,R)$:
$$ (\varphi r)(m):=\varphi(rm).$$
For all $r\in R$, $m\in M$, we have
$$\sum_irm_i\varphi^i(m)=rm=\sum_im_i\varphi^i(rm)=\sum_im_i(\varphi^ir)(m)$$
hence
\begin{equation}\label{eq1.10}\sum_irm_i\otimes\varphi^i=\sum_im_i\otimes\varphi^ir.
\end{equation}
Therefore the map $\db_{k,M}$ extends to a map $\db_M:R\to M\otimes_RM^*$
of $R$-bimodules. On the other hand, we also have an $R$-bimodule map
$\ev_M:M^*\otimes_RM\to R$, $\varphi\otimes m\mapsto \varphi(m)$,
since $\varphi(rm)=(\varphi r)(m)$. It is easy to check the following identities for 
$\ev_M$ and $\db_M$:
\begin{equation}\label{eq1.11}
(\ev_M\otimes_R\id_{M^*})(\id_{M^*}\otimes_R\db_M)=\id_{M^*};
\quad
(\id_M\otimes_R\ev_M)(\db_M\otimes_R\id_M)=\id_M.
\end{equation}
In the language of monoidal categories we call such an $R$-bimodule $M$ a
left rigid object (in $\RBimod$) and $M^*$ the left dual to $M$.

We also have the notion of right dual to a left $R$-module as well as the notion of
right rigid $R$-bimodules. In particular, the dual bimodule $M^*$ to $M$, if it exists,
is right rigid and the right dual to $M^*$ is $M$.

We define now the stucture of an $R$-coring on $M^*\otimes_kM$ for a finitely
generated projective right $R$-module $M$. Denote by $\tau$ the action of $R$
on $M^*\otimes_kM$ which is given by the action of $R$ on $M$ and denote by 
$\sigma$ the action on $M^*\otimes_kM$ which is given by the action on $M^*$.
Set
$$\Delta:=\id_{M^*}\otimes_k\db_{k,M}\otimes_k\id_{M^*}:
M^*\otimes_kM\to M^*\otimes_kM\otimes_RM^*\otimes_k M$$
and $\varepsilon:=\ev_{k,M}$. It follows immediately from \rref{eq1.9}
that $M^*\otimes_kM$ is an $R$-coring. If moreover $M$
is an $R$-bimodule then there are four actions of $R$ on $M^*\otimes_kM$.
Thus we have proved:
\begin{lem}  Let $M$ be an $R$-bimodule which is f.g. projective as a right $R$-module. 
Denote the actions of $R$ on $M^*\otimes_kM$ induced from those
on $M$ by $\tau$ and the actions of $R$ on $M^*\otimes_kM$ induced from those
on $M^*$ by $\sigma$. Then $(M^*\otimes_kM,\ev_{k,M})$ is an $R$-coalgebroid.\end{lem}

In particular, $R\otimes_kR$ is an $R$-coalgebroid, the actions of $R$ on $R\otimes_kR$
are specified as follows:
$$\si(a)\tau(b)(m\otimes n)\si(c)\tau(d)=amc\otimes bnd.$$

\subsection{Comodules over coalgebroids}\label{sec1.3}

Let $M$ be a right $R$-module and $L$ be an $R$-coalgebroid. Denote by
$\tau$ the right action of $R$ on $M$. We form the tensor product $M\tsot L$.
On this module there are three actions of $R$, induced from its actions on $L$.
A {\em coaction} of $L$ on $M$ is a map $\delta:M\lora M\tsot L$, 
satisfying the following conditions:
\begin{eqnarray*} 
\delta(m\tau(a))&=& \delta(m)\tau(a) \quad \mbox{(the linearity on $R$)},\\
(\delta \tsot \id_L)\delta&=&(\id_M \tsot \Delta)\delta \quad\mbox{(the coassociativity)},\\
(\id_M\tsot\va)\delta&=&\id_L \quad\mbox{(the unity)}.\end{eqnarray*}
In other words, $M$ is a comodule over the $R$-coring $L$ with respect
 to the ($\si,\tau$) action ($R$ acts on the left by $\si$ and on the right by $\tau$).
We use Sweedler's notation for the coaction $\delta(m)=\sum_{(m)}m_{(0)}\tsot m_{(1)}.$
Analogously, for  a left $R$-module $M$ with the action denoted by $\si$,
we can define the notion of a left coaction  of $L$ on $M$.

\example For an $R$-bimodule $M$ which is f.g. projective as right module,
$M$ is a right comodule over $L=M^*\ot_k M$ and $M^*$ is a left 
comodule over $L$.
The action is given as follows: $\delta(m)=\sum_im_i\ot_R (\varphi_i\ot_k m)$.
In particular, for $M=R$, the coaction of $R\ot_kR$ on $R$ is given
by $\delta(a)=1\ot_R (1\ot_k a)$. Note
that in this definition, we cannot move $a$ to the left, i.e., 
$\delta(a)\neq a\ot (1\ot 1)$, unless $a$ is in $k$.\bigskip

Let $L$ be an $R$-coalgebroid and $M$ a right comodule over $L$. We set
\begin{equation}\label{eq18}\tau(a)m:=\sum_{(m)}m_{(0)}\tau\va\left(\tau(a)m
_{(1)}\right).
\end{equation}
This definition does not depend on the choice of $m_{(0)}$ and $m_{(1)}$.
 Indeed, we have,
for $m\in M$, $l\in L$, $a,b\in R$,
\begin{eqnarray*} m\tau(b)\tau\va\largebr{\tau(a)l}
=m\tau\va\largebr{\si(b)\tau(a)l}
= m\tau\va\largebr{\tau(a)(\si(b)l)}.\end{eqnarray*}
\begin{lem}\label{lem1.3.2} Let $L$ be an $R$-coalgebroid.
Then the action defined in \rref{eq18}
 makes $M$ an $R$-bimodule. $\delta$ is $R$-linear with respect to
this new action on $M$ and the action on $M\tsot L$ specified above, i.e.,
$$\delta(\tau(a)m)=\sum_{(m)}m_{(0)}\ot \tau(a)m_{(1)}.$$
Furthermore, $\delta$ satisfies the equation
\begin{equation}\label{eq19} \sum_{(m)}\tau(a)m_{(0)}\ot m_{(1)}
=\sum_{(m)}m_{(0)}\ot m_{(1)}\si(a).\end{equation}
Conversely, a left action on $M$ of $R$ with respect to which $\delta$ 
is linear in the above sense
is uniquely given by the formula in \rref{eq18}.\end{lem}
 The proof contains lengthy verifications using definitions and the relations in
(\ref{eq11},\ref{eq12},\ref{eq13},\ref{eq19}) and will be omitted.

\remark  By virtue of  Lemma \ref{lem1.3.2}, by a (right) comodule over an
$R$-algebroid $L$ we shall understand
an $R$-bimodule equipped with a coaction $\delta$ satisfying 
the conditions of this lemma.
It is however not true that if $M$ is a right $L$-comodule and $N$ 
is an $R$-bimodule then
$N\ot_R M$ is an $L$-comodule for this would contradict Lemma \ref{lem1.3.2}.

Analogously, we have a notion of left $L$-comodules.
Sweedler's notation for $\delta:M\lora L\tsot M$ reads
$\delta(\varphi)=\sum_{(\varphi)}\varphi_{(-1)}\tsot \varphi_{(0)}$.
As in the case of right comodules, we can define a right action of 
$R$ on a left $L$-comodule
\begin{equation}\label{eq17.2} \varphi\si(a):=
\sum_{(\varphi)}\si\va\largebr{\varphi_{(-1)}\si(a)}\varphi_{(0)}.\end{equation}
This action is well defined and an analog of Lemma \ref{lem1.3.2} holds:
\begin{equation}\label{eq17.3}
 \sum_{(\varphi)}\tau(a)\varphi_{(-1)}\ot \varphi_{(0)}=
 \sum_{(\varphi)}\varphi_{(-1)}\ot \varphi_{(0)}\si(a).\end{equation}

\begin{lem}\label{lem1.3.3} Let $M$ be a right $L$-comodule 
which is f.g. projective
as a right module over $R$. Then there is a left coaction of $L$ 
on $M^*$ given by the condition
$$ \sum_{(\varphi)}\varphi_{(-1)}\tau\varphi_{(0)}(m)=\sum_{(m)}\si\varphi(m_{(0)})m_{(1)}.$$
This correspondence is one-to-one between right $L$-comodules,
f.g. projective as right $R$-modules and left $L$-comodules,
f.g. projective as left $R$-modules.\end{lem}
\proof Given a right $L$-comodule $M$.  The coaction on $M^*$ is given as follows:
$$\xymatrix{ M^*\ar[rr]^{\delta}\ar[d]_{\id\otimes \db_{k,M}}&&L\tsot M^*\\
M^*\ot_k M\tsot M^*
\ar[rr]_{\id\otimes\delta\otimes\id}& &M^*\ot_k(M\tsot L)\tsot M^*\ar[u]_{\ev_{k,M}\otimes\id}}$$
More explicitly, let $m_i$ and $\varphi^i$, $i=1,2,...,d$ be a pair of dual
bases in $M$ and $M^*$ respectively. The coaction on $M^*$ is given by
$$\delta(\varphi)=\sum_i\si\varphi({m_i}_{(0)})m_i{}_{(1)}\ot \varphi^i.$$
The verification is straightforward.

Conversely, given a left coaction $\delta: M\lora L\tsot M$, where $M$
is an f.g. projective left $R$-module,
one defines a right coaction of $L$ on the right dual ${}^*M$ by the condition
$\sum_{(m)} m_{(-1)}\tau\eta(m_{(0)})=\sum_{(\eta)}\si\eta_{(0)}(m)\eta_{(1)}.$
It is explicitly given as follows:
$$\xymatrix{ \rd M\ar[d]_{\db_{k,\rd M}\otimes\id} \ar[rr]^\delta&&\rd M\tsot L\\
\rd M\tsot M\ot_k \rd M\ar[rr]_{\id\otimes\delta\otimes\id}&&
 \rd M\tsot (L\tsot M)\ot_k\rd M\ar[u]_{\id\otimes\ev_{k,\rd M}} }$$
\eee
\subsection{Tensor products of coalgebroids}\label{sec1.4}
$R$-coalgebroids form a category in a natural way: morphisms between
two coalgebroids are those $R-R$-bimodules maps that commute with
$\Delta$ and $\va$. In this section, we introduce a tensor product in this
category. Let $\sqot$ denote the tensor product ${}^\si_\tau\ot^\tau_\si$,
which is given precisely by
$$L \mbox{${}^\si_\tau\ot^\tau_\si$} K=L\ot_kK\left/\Largebr{\si(a) h
 \tau(b)\ot k=h\ot \tau(b)k\si(a)}\right.,$$ for $R-R$-bimodules $L$ and $K$.
In other words, we have the following relation in $L\sqot K : \forall h\in L, k\in K$,
\begin{equation}\label{eq20} \si(a)h\tau(b)\sqot k=h\sqot \tau(b)k\si(a).\end{equation}
We specify the following actions of $R$ on $L\sqot K$:
\begin{eqnarray*} \si(a)(h\sqot k)\si(b)&=&h\si(b)\sqot \si(a) k,\\
\tau(a)(h\sqot k)\tau(b)&=&\tau(a)h\sqot k\tau(b).\end{eqnarray*}
Here we adopt the convention that the action of $R$ has preference over the tensor product.

Let $L$ and $K$ be $R$-coalgebroids. Define the $k$-linear maps
 $$\bar\Delta:L\sqot K\lora (L\sqot K) \tsot(L\sqot K)\quad\mbox{ and }\quad
 \bar \va: L\sqot K\lora R$$
as follows:
\begin{eqnarray*}
 \bar\Delta(h\sqot k)&=&\sum_{(h),(k)}\largebr{h_{(1)}\sqot k_{(1)}}
 \tsot\largebr{h_{(2)}\sqot k_{(2)}}\\
\bar\va(h\sqot k)&=&\va(k\si\va(h)) 
\end{eqnarray*}
 
The maps $\bar\Delta$ and $\bar\va$ are well defined and
 they define a coalgebroid structure on $L\sqot K$ (cf. \cite[3.10]{takeuchi2}).
Recall that $R\ot_kR$ is an $R$-coalgebroid with the $R-R$-bimodule structure
 given as follows:
$$ \si(a)\tau(b)(m\ot n)\si(c)\tau(d)=amc\ot bnd.$$
%\begin{cor}\label{cor1.3.2} 
The category of $R$-coalgebroids is a monoidal category, with the
unit object being $R\ot_k R$.

\subsection{$R$-bialgebroids} \label{sec1.5}
Since the category of $R$-coalgebroids is monoidal, we have the notion
of monoids in this category,
which are called $R$-bialgebroids. More explicitly, an $R$-bialgebroid $L$ is a coalgebroid
equipped with the following morphisms of $R-R$-bimodules
$m:L\sqot L\lora L$ and  $u:R\ot_k R\lora L$,
satisfying
\begin{eqnarray}\label{eq24}&& \Delta m=(m\tsot m)\bar\Delta; \quad \va m=\bar\va;\\
\label{eq25} &&\Delta u=u\tsot u;\quad \va u(a\ot b)=ab; \\
\label{eq26}&& m(\id_{R\ot_k R}\sqot m)=m(m\sqot \id_{R\ot_kR});\\
\label{eq27} &&m(\id_{R\ot_k R}\sqot u)= m(u\sqot \id_{R\ot_k R})=\id_L,\end{eqnarray}
where we use the identification $L\sqot (R\ot_kR)\cong L\cong (R\ot_kR)\sqot L$,
which is given explicitly by
\begin{equation} h\sqot (a\ot b)\longleftrightarrow \si(a)h\tau(b);\quad (a\ot b)\sqot
h\longleftrightarrow \tau(b)h\si(a).\end{equation}

Denoting $h\circ k=m(h\sqot k)$ and using Sweedler's notation, we have
\begin{eqnarray}\nonumber && \si(a)h\tau(b)\circ k=h\circ \tau(b)k\si(a),
\quad \va(h\circ k)=\va(k\si(h))\\
\label{eq29}&&\Delta(h\circ k)= \sum_{(h)(k)} h_{(1)}\circ k_{(1)}\tsot h_{(2)}
\circ k_{(2)};\quad \Delta(1_L)= 1_L\tsot 1_L, \end{eqnarray}
bearing in mind the preference of $\circ$ over $\ot$, where $1_L:=u(1_R\ot_k 1_R)$.
On the other hand, the linearity of $m$ and $u$ can be expressed as
\begin{equation}\label{eq30}\si(a)\tau(b)(h\circ k)\tau(c)\si(d)=
 \tau(b)h\si(d)\circ \si(a)k\tau(c),\end{equation}
where  we adopt the convention that the action of $R$ has preference over the product.

We define maps $s$ and $t$ from $R$ to $L$ as follows
$$ s(a)=\si(a)1_L=1_L\si(a); \quad t(a)=\tau(a)1_L=1_L \tau(a).$$
The following relations follow immediately from \rref{eq27} (or from \rref{eq30})
\begin{eqnarray} \label{eq31}\begin{array}{cc} s(a)\circ h=h\si(a); &h\circ s(a)=\si(a)h;\\
 h\circ t(a)=h\tau(a);& t(a)\circ h=\tau(a)h.\end{array}\end{eqnarray}
In particular, $s$ is an anti-homomorphism and $t$ is a homomorphism of $k$-algebras from $R\lora L$.

Bialgebroids over associative algebras seem to be first introduced by Takeuchi \cite{takeuchi1} and later
independently introduced by J. Lu \cite{lu1}.
\subsection{Comodules over bialgebroids}\label{sec1.6}
A comodule over a bialgebroid is by definition a comodule over the underlying coalgebroid. We
have seen in Subsection \ref{sec1.3} that a right comodule over an $R$-coalgebroid, which is initially
a right $R$-module, can be endowed with a structure of left $R$-module.
In this subsection we show that the tensor product of two comodules
over a bialgebroid is again a comodule.

Let $M,N$ be right comodules over an $R$-bialgebroid $L$.
Define a coaction of $L$ on $M\ot_R N$ as follows:
$$ \delta(m\ot n)=\sum_{(m)(n)}m_{(0)}\ot n_{(0)}\ot m_{(1)}\circ n_{(1)}.$$
\begin{lem}\label{lem1.4.1} The coaction given above is well defined and makes
$M\ot_R N$ a comodule over $L$.\end{lem}

Notice that $R$ itself is a comodule over $L$ by means of the morphism $t$
defined in Subsection \ref{sec1.5}:
$\delta(a)=1\ot t(a)=1\ot \tau(a)1_H$.
\begin{cor}\label{cor1.6.2} The category of comodules
over a bialgebroid is monoidal with the unit object being $R$.\end{cor}

\subsection{The antipode}\label{sec1.7}
Consider the tensor product $H\ssot H$ defined as follows
$$ H\ssot H:= H\ot_k H\left/\Largebr{\si(a)h\ot k=h\ot k\si(a)}\right.$$
and specify the actions of $R$ as follows
$$ \tau(a)(h\ot k)\tau(b)=\tau(a)h\ot k\tau(b);\quad \si(a)(h\ot k)\si(b)=h\si(b)\ot \si(a)k.$$
There is an $R-R$-bimodule morphism $\pi:H\ssot H\lora H\sqot H,$ which is a quotient map.

\definition (cf. \cite{schauen2}) \rm  Let $H$ be an $R$-bialgebroid. An antipode on $H$ is by
definition a map $\nabla:H\lora H\ssot H; \nabla(h)=\sum_{(h)}h^-\ssot h^+,$
satisfying the following conditions:
\begin{eqnarray}\label{eq341}\sum_{(h)}h^-\circ {h^+}_{(1)}\tsot {h^+}_{(2)}&=&1\tsot h\\
\label{eq342}\sum_{(h)} h_{(1)}\circ {h_{(2)}}^-\ssot{h_{(2)}}^+&=&1\ssot h.\end{eqnarray}
If such an antipode exists, $H$ is called a Hopf algebroid.\\

Define a map $\beta:H\ssot H\lora H\tsot H$ to be
\begin{equation} \label{eq34}\xymatrix{ 
H\ssot H\ar[d]_{\id\ot \Delta} \ar[rr]^\beta&&H\tsot H\\
H\ssot (H\tsot H)\ar[r]^\cong & (H\ssot H)\tsot H\ar[r]^\pi
&
(H\sqot H)\tsot H\ar[u]_{m\otimes\id} 
 }\end{equation}
$$\beta\left(h\ssot k\right)=\sum_{(k)}h\circ k_{(1)}\tsot k_{(2)}.$$
Then we have
$$\textstyle \beta\left(\sum_{(k)}h\circ k^-\ssot k^+\right)=
\sum_{(k)}h\circ k^-\circ {k^+}_{(1)}\tsot {k^+}_{(2)}=h\tsot k
\mbox{ by \rref{eq341} }$$
and
$$ \sum_{(k)}h\circ k_{(1)}\circ{k_{(2)}}^-\ssot {k_{(2)}}^+=h\tsot k\mbox{ by \rref{eq342} }.$$
Therefore the map $\beta$ is invertible with the inverse given by
\begin{equation}\label{eq343}\beta^{-1}(h\tsot k)=\sum_{(k)}h\circ k^-\ssot k^+.\end{equation}
We have
$ \nabla(h)=\beta^{-1}(1\ot h),$
whence $\nabla$ is uniquely determined. Thus, if an antipode exists then it is determined uniquely.

\remark
 If $R=k$ and $H$ is a Hopf algebra over $k$ 
then $\nabla$ is given explicitly by $\nabla(h)=\sum_{(h)}S(h_{(1)})\ot h_{(2)}$
where $S$ denotes the antipode of $H$.

\begin{lem}[{\cite[Pro.~3.7]{schauen2}}] \label{lem7.1}Let $H$ be a Hopf algebroid.
Then the antipode $\nabla$
satisfies the following relations:
\begin{eqnarray}\label{eq35.0}&&\nabla(1_H)=1_H\ot 1_H,\\
\label{eq35} &&\nabla(\tau(a)\si(b)h\si(c)\tau(d))
=\sum_{(h)}\tau(b)h^-\tau(c)\ot \tau(a)h^+\tau(d),\\
 \label{eq39}&& \sum_{(h)}{h_{(1)}}^-\ssot {h_{(1)}}^+\tsot h_{(2)}=
\sum_{(h)}h^-\ssot {h^+}_{(1)}\tsot h{^+}_{(2)},\\
&&\label{eq41}\sum_{(h)}h^-\stdot h^{+-}\ssot h^{++}=
\sum_{(h)}h^-{}_{(2)}\stdot h^-{}_{(1)}\ssot h^+,\\
&&\label{eq41b} \sum_{(h)} h^+\si(\va(h^-))=h,\\
&&\label{eq41c} \sum_{(h)}h^-\circ h^+=1_H\tau\va(h).\end{eqnarray}
\end{lem}
\proof Applying $\beta$ on both sides of Eq. \rref{eq35.0} and \rref{eq35}
we obtain the identity maps. Thus, the equalities follows from the invertibility
of $\beta$.

Applying $\beta$ on the first two tensor components of both sides of
Eq. \rref{eq39} and using \rref{eq341}, we obtain the same values. 
Thus the equality also follows
from the invertibility of $\beta.$ 

We prove \rref{eq41}. Let $\bar\beta$ be the map
$H\stdot H\ssot H\lora H\tsot H\tsot H$; 
$$\bar\beta(h\stdot k\ssot l)=k\circ l_{(1)}\tsot h\circ l_{(2)}\tsot l_{(3)}$$
$\bar\beta$ is also invertible with the inverse given by 
$\bar\beta^{-1}(k\tsot h\tsot l)=h\circ l^-\stdot k\circ l^{+-}\ssot l^{++}.$
Applying $\bar\beta$ to the right-hand side of \rref{eq41} we obtain:
\begin{eqnarray*}\lefteqn{\sum_{(h)} h^-{}_{(1)}\circ  h^+{}_{(1)} \tsot 
 h^-{}_{(2)}\circ
h^+{}_{(2)}\tsot h^+{}_{(3)}}\\ 
&=& 
(h^-\circ h^+{}_{(1)})_{(1)}\tsot (h^-\circ h^+{}_{(1)})_{(2)}\tsot h^+{}_{(2)}
\quad\mbox{by \rref{eq29}}\\
&=& 1\tsot 1\tsot h\end{eqnarray*}
On the other hand, applying $\bar\beta$ to the left-hand side of 
\rref{eq41}, we obtain
\begin{eqnarray*} 
\lefteqn{\sum_{(h)}(h^+)^-\circ (h^+)^+{}_{(1)}\tsot h^-\circ (h^+)^+{}_{(2)}\tsot (h^+)^+{}_{(3)}
}\\ &=&
\sum_{(h)} 1\tsot h^-\circ h^+{}_{(1)}\tsot h^+{}_{(2)}\quad \mbox{by \rref{eq341}}\\
&=& 1\tsot 1\tsot h.\end{eqnarray*}
The invertibility of $\bar\beta$ implies the equality in \rref{eq41}.

Eq. \rref{eq41b} follows from Eq. \rref{eq342}. Indeed, we have
\begin{eqnarray*} h&=& \sum_{(h)}h_{(2)}{}^+\si\va(h_{(1)}\circ h_{(2)}{}^-)\\
&=& h_{(2)}{}^+\si\va\largebr{h_{(2)}{}^-\si\va(h_{(1)})}\\
&=& h_{(2)}{^+\si\va\largebr{\tau\va(h_{(1)}h_{(2)}{}^-)}} 
\quad \mbox{ by condition (iii) for $\va$}\\
&=& h^+\si\va(h^-) \quad \mbox{ by \rref{eq35} for $\tau(b)$ and by \rref{eq10}}
\end{eqnarray*}

Eq. \rref{eq41c} also follows immediately from \rref{eq341} and \rref{eq10}.
\eee

\begin{pro}\label{pro7.1} Let $H$ be an $R$-Hopf algebroid and $M$ a right $H$-comodule.
 Assume that $M$ is an f.g. projective  right $R$-module. Then there exists a coaction of $H$
on $M^*$ making it the dual object to $M$ in the category of right $H$-comodules.\end{pro}
\proof
We first define a coaction of $H$ on $M^*$. The right coaction of $H$ on $M$
induces the left coaction of $H$ on $M^*$: 
$\partial (\varphi)=\sum_{(\varphi)}\varphi_{(-1)}\tsot \varphi_{(0)}$
by the condition
\begin{equation}\label{eq43.0} \sum_{(\varphi)}\varphi_{(-1)}\tau \varphi_{(0)}(m)=
\sum_{(m)}\si\varphi(m_{(0)})m_{(1)}.\end{equation}
Define a right coaction of $H$ on $M^*$ as follows
\begin{equation}\label{eq43} \delta(\varphi):=
\sum_{(\varphi)}\si\va\largebr{\varphi_{(-1)}{}^+}\varphi_{(0)}\ssuot \varphi_{(-1)}{}^-.\end{equation}
We will show that this is a well defined coaction of $H$ and that with this
coaction $M^*$ is a dual $H$-comodule
to $M$.

It is easy to check that this coaction  is well defined, i.e.,
it does not depend on the choice of $\varphi_{(0)}, \varphi_{(-1)}$ and
$\varphi_{(-1)}{}^-, \varphi_{(-1)}{}^+$.
We show that $\delta$ is a coaction. For simplicity we shall use the 
notation $ \phi:=\varphi_{(-1)}$. Notice that
$$\sum_{(\varphi)}\varphi_{(-2)}\tsot \varphi_{(-1)}\tsot\varphi_{(0)}=
\sum_{(\varphi)(\phi)}\phi_{(1)}\tsot \phi_{(2)}\tsot \varphi_{(0)}.$$

The coassociativity of $\delta$ amounts to the following equation
\begin{eqnarray*}
\lefteqn{\sum_{(\varphi)}\si\va\largebr{\phi_{(2)} {}^+}\varphi_{(0)}\ssuot
\tau\va\largebr{\phi_{(1)}{}^+}\phi_{(2)} {}^-\tsot\phi_{(1)}{}^-}
\qquad\qquad && \\
&=& \sum_{(\varphi)}\si\va\largebr{\phi {}^+}\varphi_{(0)}\ssuot
\largebr{\phi {}^-}_{(1)}\tsot \largebr{\phi {}^-}_{(2)}
\end{eqnarray*}
Applying $\nabla$ on the last term of \rref{eq39} 
and taking in account \rref{eq35}, we obtain the equality 
$$\sum_{(h)} h_{(1)}{}^-\ssot h_{(1)}{}^+\ttot  h_{(2)}{}^-\ssot h_{(2)}{}^+=
\sum_{(h)}h^-\ssot h^+{}_{(1)}\ttot (h^+{}_{(2)})^-\ssot (h^+{}_{(2)})^+.$$
 Therefore, for $h=\phi=\varphi_{(-1)}$, we have
\begin{eqnarray*}
\lefteqn{\sum_{(\varphi)}\si\va\largebr{\phi_{(2)} {}^+}\varphi_{(0)}\ssuot
\tau\va\largebr{\phi_{(1)}{}^+}\phi_{(2)} {}^-\tsot\phi_{(1)}{}^-}&& \\
&=&\sum_{(\varphi)}\si\va\Largebr{\largebr{\phi {}^+{}_{(2)}}^+} \varphi_{(0)}\ssuot
\tau\va\Largebr{{\phi {}^+}_{(1)}} \largebr{\phi {}^+{}_{(2)}}^-\tsot \phi {}^-\\
&=&\sum_{(\varphi)}\si\va\Largebr{\largebr{\si\va\Largebr{{\phi {}^+}_{(1)}}\phi {}^+{}_{(2)}}^+} \varphi_{(0)}\ssuot
 \largebr{\si\va\Largebr{{\phi {}^+}_{(1)}}\phi {}^+{}_{(2)}}^-\tsot \phi {}^-\quad \mbox{ (by \rref{eq35})}\\
&=& \sum_{(\varphi)}\si\va\largebr{(\phi^+) {}^{+}}\varphi_{(0)}\ssuot
(\phi^+) ^{-}\tsot \phi {}^-\quad\mbox{(by \rref{eq10})}\\
&=& \sum_{(\varphi)}\si\va{\largebr{\phi {}^+}}\varphi_{(0)}\ssuot
\largebr{\phi {}^{-}}_{(1)}\tsot \largebr{\phi {}^{-}}_{(2)} \quad\mbox{(by \rref{eq41})}
\end{eqnarray*}

The counity amounts to the following equation
$$ \varphi=\sum_{(\varphi)}\si\va(\phi {}^+)\varphi_{(0)}\si\va(\phi {}^-).$$
Indeed,  by \rref{eq17.2}
\begin{eqnarray*}\lefteqn{ \sum_{(\varphi)}\Largebr{\si\va\largebr{\phi {}^+}\varphi_{(0)}}
\si\va(\phi {}^-)}\qquad\qquad && \\
&=&\sum_{(\varphi)} \si\va\left[\Largebr{    \si\va\largebr{\phi_{(1)}{}^+}\phi_{(2)} }
 \si\va\largebr{\phi_{(1)}{}^-}\right]    \varphi_{(0)}    \\
&=&\sum_{(\varphi)} \si\va\left[\Largebr{    \si\va\largebr{\phi {}^+{}_{(1)}}
\phi {}^+{}_{(2)}} \si\va\largebr{\phi {}^-}\right]    \varphi_{(0)}  \quad\mbox{by \rref{eq39}}  \\
&=&\sum_{(\varphi)} \si\va\largebr{    \phi {}^+\si\va\largebr{\phi {}^-}}
 \varphi_{(0)} \quad\mbox{by \rref{eq41b}}   \\
&=& \sum_{(\varphi)}\si\va\largebr{\varphi_{(-1}}\varphi_{(0)}=\varphi
\end{eqnarray*}
Similar computation shows that the left action of $R$ on $M^*$, which is induced from
 the right coaction of $H$ as in \rref{eq18} is just the natural one: 
 $(a,\varphi)\loma \si(a)\varphi: [\si(a)\varphi](m)=a\varphi(m)$. Thus, Lemma \ref{lem1.3.2} 
 applies and the equation in \rref{eq19} has the form
 \begin{equation}\label{eq43b}\si(a)\si\va(\phi^+)\varphi_{(0)}\ssuot  \phi^-=
 \si\va(\phi^+)\varphi_{(0)}\ssuot \phi^-\si(a)
 \end{equation}
 
Finally, we show that $M^*$ equipped with this coaction is
 a left dual comodule to $M$, which amounts
to showing that $\ev:M^*\stot M\lora R$ and $\db:R\lora M\tsot M^*$
 are morphisms of $H$-comodules. 
Choose a pair of dual bases on $M$ and $M^*$,
 $\{m_i\}$ and $\{\varphi^i\}$, and
denote for simplicity $\phi:=\varphi_{(-1)}$ and $\phi^i:=\varphi^i{}_{(-1)}$.
We have to check the following equations:
\begin{eqnarray}\label{eq45}
\sum_{(\varphi)} \va\largebr{\phi {}^+}\varphi_{(0)}(m_{(0)})\tsot
\largebr{\phi {}^-\circ m_{(1)}}=1\tsot 1\tau\varphi(m),\end{eqnarray}
\begin{equation}\label{eq46}
\sum_{i,(\varphi^i)} m_{i(0)}\tsot \si\va\largebr{\phi^i {}^+}\varphi^i{}_{(0)}
\tsot\largebr{m_{i(1)}\circ \phi^i {}^-}
=\sum_im_i\tsot \varphi^i\tsot 1_H\end{equation}
Notice that \rref{eq46} is equivalent
 to the following: for all $\varphi\in M^*$,
\begin{equation}\label{eq47} \sum_{i,(\varphi^i)} \si\varphi(m_{i(0)})\si\va
\largebr{\phi^i {}^+}\varphi^i{}_{(0)} \tsot \largebr{m_{i(1)}\circ\phi^i {}^-}=
\varphi\tsot 1_H.\end{equation}

We prove \rref{eq45}:
\begin{eqnarray*} \lefteqn{
\sum_{(\varphi)}\va\largebr{\phi {}^+}\varphi_{(0)}(m_{(0)})\tsot \largebr{\phi {}^-\circ m_{(1)}}}\\
&=& \sum_{(\varphi)}\va\largebr{\phi {}^+}
\tsot \Largebr{\phi {}^-\circ \si\varphi_{(0)}(m)m_{(1)}}
\quad\mbox{ by \rref{eq30}}\\
&=& \sum_{(\varphi)}\va\largebr{\phi_{(1)}{}^+}\tsot
 \Largebr{\phi_{(1)}{}^-\circ \phi_{(2)} \tau\varphi_{(0)}(m)}\quad\mbox{by \rref{eq43.0}}\\
&=& \sum_{(\varphi)}\va\largebr{\phi {}^+{}_{(1)}}\tsot
 \Largebr{\phi {}^-\circ \phi {}^+{}_{(2)}\tau\varphi_{(0)}(m)}\quad\mbox{by \rref{eq39}}\\
&=&\sum_{(\varphi)}1\tsot \Largebr{\phi {}^-\circ\Largebr{\si\va\largebr{\phi {}^+{}_{(1)}}
\phi {}^+{}_{(2)} }}
\tau\varphi_{(0)}(m)\quad\mbox{by \rref{eq30}}\\
&=& \sum_{(\varphi)}1\tsot (\phi^-\circ \phi^+)\tau\varphi_{0}(m) \quad\mbox{by \rref{eq41c}}\\
&=&1\tsot 1\tau\varphi(m)
\end{eqnarray*}
For \rref{eq47}, we first notice that, on applying $\delta$ on both sides of the equations
$\varphi=\sum_i \si\varphi(m_i)\varphi^i$, we have
\begin{equation}\label{eq48} \sum_{(\varphi)}\si\va\largebr{\phi {}^+}\varphi_{(0)}\ot \phi {}^-=
\sum_{i,(\varphi^i)}\si\va\largebr{\phi^i{}^+}\varphi^i{}_{(0)}\ot
\tau\varphi(m_i)\phi^i {}^-.\end{equation}
Now, the left hand side of \rref{eq47} is equal to
\begin{eqnarray*}
\lefteqn{\sum_{i,(\varphi)} \si\va\largebr{\phi^i {}^+}\varphi^i{}_{(0)}\tsot
\Largebr{m_{i(1)}\circ \phi^i {}^-\si\varphi(m_{i(0)})} \quad\mbox{by \rref{eq43b}
%\rref{eq17.3}, \rref{eq43}
}}&&\\
&=&\sum_{i,(\varphi)} \si\va\largebr{\phi^i {}^+}\varphi^i{}_{(0)}\tsot
\Largebr{\si\varphi(m_{i(0)})m_{i(1)}\circ \phi^i {}^-}\quad \mbox{by \rref{eq29}}\\
&=& \sum_{i,(\varphi)}\si\va\largebr{\phi^i {}^+}\varphi^i{}_{(0)}\tsot
 \Largebr{\phi \tau\varphi_{(0)}(m_i)\circ \phi^i {}^-}\quad\mbox{by \rref{eq43.0}}\\
&=&\sum_{i,(\varphi)} \si\va\largebr{\phi^i {}^+}\varphi^i{}_{(0)}\tsot
\Largebr{\phi \circ\tau\varphi_{(0)}(m_i)\phi^i {}^-}\\
&=&\sum_{i,(\varphi)}\si\va\largebr{\phi^i {}^+}\varphi^i{}_{(0)}\tsot
 \Largebr{\phi_{(1)}\circ \phi ^-}\quad\mbox{by \rref{eq48}}\\
&=&\sum_{(\varphi)} \si\va\largebr{\phi }\varphi_{(0)}\tsot 1_H\quad\mbox{by \rref{eq48}}\\
&=& \varphi\tsot 1_H\end{eqnarray*}
The proof is complete.\eee
\subsection{The opposite antipode}
We have seen that for a Hopf algebroid, each comodule which is f.g. projective over $R$
possesses a left dual. As we know in the case of Hopf algebras, a right dual can be defined in terms
of the inverse to the antipode, i.e., if the antipode is bijective,
each finite dimensional comodule possesses
a right dual. A (sufficient) condition for the existence of the right dual to a f.g. projective comodule
of a Hopf algebroid  can be expressed as the existence of
a map generalizing the map $h\loma S^{-1}(h_{(2)})\ot h_{(1)}$ for Hopf algebras 
(see Remark in Subsection \ref{sec1.7}).

\definition \rm Let $H$ be a bialgebroid. An {\em opposite antipode} is a map
$\nablaop:H\lora H\ttot H$, $\nablaop(h):=\sum_{(h)} h_-\ttot h_+$, satisfying the following axioms:
\begin{eqnarray*}
 \sum_{(h)} h_{+(1)}\tsot h_-\circ h_{+(2)}=h\tsot 1\\
 \sum_{(h)} h_{(2)}\circ h_{(1)-}\ttot h_{(1)+}=1\ttot h.\end{eqnarray*}

\begin{lem}\label{lem7.2}
Let $H$ be a bialgebroid with an  opposite antipode $\nablaop$. Define a map
$\gamma:H\ttot H\lora H\tsot H$,
$\gamma(h\ot k)=\sum k_{(1)}\tsot h\circ k_{(2)}$. Then $\gamma$ is invertible with the inverse given
by $\gamma^{-1}(h\tsot k)= \sum k\circ h_-\ttot h_+$. Further the map $\nablaop$ satisfies the following
equations:
\begin{eqnarray*}&& \nablaop(\tau(a)\si(b)h\si(c)\tau(d)=\sum_{(h)}\si(a)h_-\si(d)\ttot \si(b)h_+\si(c)\\
&&\sum_{(h)} h_{+(1)}\tsot h_{+(2)}\ttdot h_-=\sum_{(h)} h_{(1)}\tsot h_{(2)+}\ttdot h_{(2)-}\\
&&\sum_{(h)} h_+\tsot h_{+-}\ttot h_{++}=\sum_{(h)} h_{-(1)}\tsot h_{-(2)}\ttot h_+.\end{eqnarray*}
Consequently, the opposite antipode is determined uniquely.
\end{lem}
We call a bialgebroid equipped with an opposite antipode {\em opposite Hopf algebroid}.
\begin{pro}\label{pro7.2} Let $H$ be an opposite $R$-Hopf algebroid and $\delta: M\lora M\tsot H$ a right
coaction of $H$ on $M$. Then the opposite antipode induces a left coaction of $H$
on $M$, $M\lora H\ttot M$,
given by
$$ m\loma \sum m_{(1)-}\ttot m_{(0)}\tau\va(m_{(1)+}).$$
Assume that $M$ is an f.g. projective left $R$-module with the right dual $\rd M$. Then
there exists a coaction $\rho:\rd M\lora \rd M\ssuot H$, $\eta\loma \sum_{(\eta)} \eta_{(0)}\ot \eta_{(1)}$,
 making it a right dual $H$-comodule to $M$.
$\rho$ is given by the following condition
$$\sum_{(m)(\eta)} m_{(1)-}\tau\eta(m_{(0)}\tau\va(m_{(1)+}))=
\sum_{(m)(\eta)}\si\eta_{(0)}\largebr{m_{(0)}\tau\va(m_{(1)})}\eta_{(1)}.
$$
\end{pro}
The proof of these facts is left to the reader.
\section{The Tannaka-Krein duality}\label{sec2}
\subsection{Tannaka-Krein duality for corings}\label{sec2.1}
We fix as in Section 1 a commutative ring $k$ and assume that everything
is $k$-linear. Let $R$ be a $k$-algebra. Tannaka-Krein duality for
 $R$-corings was proved by P. Deligne
\cite{deligne90}. Our presentation here follows A. Brugui\`eres \cite{brug1}.

Let $\cat C$ be a category and $\fun F:\cat C\lora \ModR$ be a functor
to the category of right $R$-modules.  We define the $\coend$ of $\fun F$
to be an $R$-bimodule $L$ satisfying the following universal property:
for any  $R$-bimodule $C$, there is a natural isomorphism
\begin{equation}\label{rt1} \Nat_{R}(\fun F,\fun F\ot_R C)\cong {}_R\Hom_{R}(L,C).\end{equation}
 Here we use the convention of Subsection \ref{sec1.1} for the $\Hom$.
 By the universal property, $L$, if it exists, is uniquely determined up to an isomorphism.

If the image of $\fun F$ lies in the subcategory of right f.g. projective $R$-modules,
$L$ can be constructed as follows. First notice that for any $R$-bimodule $C$ and any
object $X\in\cat C$, we have by means of the projectivity of $\fun F(X)$
$$\Hom_R(\fun F(X),\fun F(X)\otimes_RC)\cong {}_R\Hom_R(
\fun F(X)^*\otimes_k\fun F(X),C)$$
Thus we form the direct sum
\begin{equation}\label{rt2} L_0:=\bigoplus_{X\in\cat C}
\fun F(X)^*\otimes_k\fun F(X)
\end{equation}
and for any morphism $f:X\to Y$ in $\cat C$, consider the (inner) diagram of $R$-bimodule maps:
\begin{equation}\label{rt3}
\xymatrix{\fun F(Y)^*\otimes_k\fun F(X)\ar[r]^{\fun F(f)^*\otimes\id}
\ar[d]^{\id\otimes\fun F(f)}&
\fun F(X)^*\otimes_k\fun F(X)\ar[d]\ar[rdd]\\
\fun F(Y)^*\otimes_k\fun F(Y)\ar[r]\ar[drr]&L_0\ar@{-->}[rd]\\
&&L}
\end{equation}
Let $L$ be the maximal quotient $R$-bimodule of $L_0$ which makes all
the above (outer) diagrams commute. Then it is easy to see that the $R$-bimodule $L$ satisfies the universal property
in \rref{rt1}. 

 We will show that $L$ is an $R$-coring in the sense
of Subsection \ref{sec1.1}. As in  \ref{sec_example} we denote the
actions of $R$ on $\fun F(X)$
by $\tau$ and the actions on its dual by $\si$. Thus
$\fun F(X)^*\otimes_k\fun F(X)$ is an $R$-bimodule by means of the actions $\tau$ and $\si$.
The actions of $R$ on $L_0$ and $L$ will be named in the same way.

Set $C=L$ in \rref{rt1}. Then the identity $L\to L$ corresponds though
the isormophism in \rref{rt1} to  a natural transformation
$\delta:\fun F\lora \fun F\ot_RL$. For an arbitrary natural transformation
$\rho:\fun F\lora \fun F\ot_RC$, the naturality of \rref{rt1} on $C$
implies that the corresponding morphism $f_\rho:L\lora C$ satisfies
 \begin{equation}\label{rt4}\rho=(\id\ot f_\rho)\delta.
\end{equation}
 For $C=L\tsot L$ the natural morphism 
 $$(\delta\ot \id)\delta:\fun F\to\fun F\otimes_R L\tsot L$$
 corresponds though the isomorphism in \rref{rt1} to 
  a morphism $\Delta:L\lora L\tsot L$, which according to \rref{rt4} satisfies
  $$(\delta\ot \id)\delta=(\id\ot \Delta)\delta.$$
  For $C=R$, the identity transformation corresponds to a morphism: $\va:L\lora R$.
It is easy to show that $(L,\Delta,\va)$ is an $R$-coring.

\begin{lem}\label{lem_rt2} Through the isomorphism in \rref{rt1}, if $\delta\in
\Nat(\fun F,\fun F\ot_R C)$ is a family of coactions of an $R$-coring $C$, then
the corresponding morphism $L\lora C$  is a morphism of $R$-corings.\end{lem}
\proof  Let $\varphi:L\lora C$ be the map that corresponds to $\delta$.
Since $\delta$ is a coaction of $C$ on $\fun F(X)$ for every $X$, we have two
equal maps $(\delta\ot\id)\delta=(\id\ot\Delta)\delta:\fun F(X)\lora \fun F(X)\ot C\ot C$.
These maps should correspond to the same map $L\lora C\ot C$, which, that is
$\Delta_C\varphi=(\varphi\ot\varphi)\Delta_L$. The commutativity of $\varphi$ with the counits also
follows from the universal property of $L$.\eee

Thus, given a category $\cat C$ and a  functor
$\fun F:\cat C\lora \ModR$ with image in the subcategory
of f.g. projective modules, then $\fun
F$ factors through a functor $\bar{\fun F}:\cat C\lora \comodL$,
 and the forgetful functor. This is the
first part of Tannaka-Krein duality. The second part, which is
usually more difficult, is to prove that $\bar{\fun F}$ is an equivalence
if $\cat C$ is a ``good'' abelian category.

\bigskip

From now on we shall assume that $k$ is a field.
Recall that a $k$-linear abelian category $\cat C$ is  said to be \emph{locally finite}
(over $k$) if each its Hom-set is finite dimension over $k$ and each object has the
composition series of finite length.

\definition Let $R$ be a $k$-algebra and $L$ be an
$R$-coring. $L$ is said to be (right) semi-transitive if the following conditions
are satisfied:
\begin{enumerate}\item each right $L$-comodule is projective over $R$.
\item each $L$-comodule is a filtered limits of subcomodules which are
finitely generated over $R$.
\item the category $\comod$-$L$ of right $L$-comodules which are finitely
generated as  $R$-modules is locally finite over $k$.\end{enumerate}

\begin{thm}[\cite{deligne90}, see also {\cite[Thm. 5.2]{brug1}}] \label{thm_deligne}
Let $k$ be a field and $\cat C$ be
a (small) $k$-linear abelian category which is locally finite. Let  $\fun F:\cat C
\lora \modR $ be an exact faithful functor with image in the subcategory
of f.g. projective modules. Let $L=\coend(\fun F)$. Then $L$ is a semi-transitive
coring and the functor $\bar{\fun F}$ is an equivalence of
abelian categories. Conversely, let $L$ be a semi-transitive coring and
$\fun F:\comodL\lora \modR $ be the forgetful functor. Then
$\fun F$ is faithful, exact and  has image in the category of projective modules of finite
rank and $L\cong \coend(\fun F)$.\end{thm}

\subsection{Tannaka-Krein duality for bialgebroids}\label{sec2.2}
Let $\cat C$ be a $k$-linear category and $\fun F:\cat C\lora \RBimod$ a functor
with image in the subcategory of left rigid $R$-bimodules (i.e., f.g. projective as
right $R$-modules). Then we can construct the $\coend$ of $\fun F$, denoted by $L$.
There are several actions of $R$ on $L$ which we will specify now.

Recall from Subsection \ref{sec_example} that the left dual $\fun F(X)^*$ to $\fun F(X)$ is also
an $R$-bimodule. We shall use the convention of  \ref{sec_example} for denoting
the actions of $R$ on $\fun F(X)^*\otimes_k\fun F(X)$. The actions of $R$ on
$L_0$ will be denoted accordingly. Since the maps $\fun F(f)^*\otimes \id$
and $\id\otimes \fun F(f)$ in the diagram \rref{rt3} commute with all the (left and right)
actions $\sigma$ and $\tau$, there are natural actions of $R$ on $L$ which will be
denoted accordingly. As shown in Subsection \ref{sec2.1}, $L$ with respect to the bimodule
structure  given by $(\sigma,\tau)$ is an $R$-coring.

\begin{lem}\label{lem_coalgebroid} Let $\fun F:\cat C\lora \RBimod$ be a functor
with image in the subcategory of left rigid bimodules. Then $L=\coend(\fun F)$ is
a coalgebroid.\end{lem}
\proof As shown in the previous subsection, $L$ is an $R$-coring with respect
to the actions $(\si,\tau)$. It remains to show that $\Delta$ is a morphism of
$R-R$-bimodules and that $\var$ satisfies $\va(\tau(r)a)=\var(a\si(r)),\  \forall r\in R,
a\in L$.

To see that $\Delta$ is a morphism of $R-R$-bimodules, it is sufficient to notice that
in the construction of $L$ (diagrams in \rref{rt3}) all maps are $R-R$-bimodules
morphisms and for each object $X\in \cat C$, the coproduct
$$\Delta_X:\fun F(X)^*\ot_k\fun F(X)\lora
\fun F(X)^*\ot_k\fun F(X)\ot_R\fun F(X)^*\ot_k\fun F(X)$$
is a morphism of $R-R$-bimodules, for $\fun F(X)^*\ot_k\fun F(X)$ is a
coalgebroid (see \ref{sec_example}).

Similarly, the counit $\var_X:\fun F(X)^*\ot_k\fun F(X)\lora R$ satisfies
$\va_X(\tau(r)a)=\var_X(a\si(r))$ and moreover, for any pair of objects $X,Y\in\cat C$, a morphism $f:X
\lora Y$ induces a morphism $\var_f:\fun F(Y)^*\ot_k\fun F(X)\lora R$ which is linear
with respect to the actions $(\si,\tau)$ and satisfies $\var_f(\tau(r)a)=\var_f(a\si(r))$.
Therefore we have commutative diagrams of the form
\begin{equation}\label{rt3.1}\xymatrix{
\fun F(Y)^*\ot_k \fun F(X)\ar[rr]^{\fun F(f)\ot \id}\ar[rd]^{\var_f}\ar[dd]^{\id\ot\fun F(f)}&&
\fun F(X)^*\ot_k\fun F(X)\ar[ld]_{\var_X} \ar[dd]\\
&R&\\
\fun F(Y)^*\ot_k \fun F(Y)\ar[ru]^{\var_Y}\ar[rr]&&L\ar@{-->}[lu]^{\var}}\end{equation}
By construction, $L$ is a quotient of $L_0$, which is the direct sum of $\fun F(X)^*\otimes_k\fun F(X)$,
$X\in \cat C$. We therefore conclude that the induced map $\var:L\to R$ also satisfies
the equation $\var(\tau(r)a)=\var(a\si(r))$.  Thus, $L$ is an $R$-coalgebroid. \eee
\begin{lem}\label{lem_rt2.1}
Let $C$ be an $R$-coalgebroid. Let $\delta\in \Nat(\fun F,\fun F\ot_RC)$
be a natural transformation, which is a family of
coactions of a coalgebroid $C$ satisfying equation \rref{eq19}. Then $\delta$
corresponds though the isomorphism in \rref{rt1}
to a morphism $L\lora C$ of coalgebroids.\end{lem}
\proof A coaction of a coalgebroid $C$ on a left rigid bimodule
$\fun F(X)$, which satisfies the equation \rref{eq19}, induces a  morphism
of coalgebroids $\bar\delta:\fun F(X)^*\ot_k\fun F(X)\lora C$. Consequently the map
$$\sum_{X\in\cat C}\bar\delta_X:L_0=\bigoplus_{x\in\cat C}\fun F(X)^*\otimes_k\fun F(X)
\lora C$$
is also a homomorphism of $R$-coalgebroids. On the other hand, 
according to Lemma \ref{lem_rt2}, there exists a homomorphism of $R$-corings
$\varphi:L\lora C$ which fits in the following commutative diagrams for
all $X\in \cat C$:
$$\xymatrix{
L_0\ar[r]\ar[rd]_{\sum_{X\in\cat C}\bar\delta_X} &L\ar[d]^\varphi\\
& C } $$
Since the map $L_0\lora L$  is surjective, the $R$-linearity (with respect to all actions)
of $\varphi$ follows from the $R$-linearity of the maps $L_0\lora L$ and $L_0\lora C$.
Thus $\varphi$ is a homomorphism of $R$-coalgebroids.
\eee

\bigskip

We recall that the tensor product $\sqot$ was introduced in \ref{sec1.4}.
\begin{pro}\label{pro_rt3} Let $\fun F$ and $\fun G$ be functors $\cat C\lora
\RBimod$ with images in the category of left rigid bimodules.
Let $L=\coend( \fun F)$ and $K=\coend( \fun G)$. Then
\begin{equation}\label{rt51} \coend (\fun F\ot_R \fun G)\cong L\sqot K.\end{equation}\end{pro}
\proof  We still keep the notation for the actions of $R$ on $M^*\otimes_kM$, $M\in\RRbimod$,
as in Subsection \ref{sec_example}. We notice the following isomorphism for the $\sqot$-product
\begin{eqnarray}\label{rt6} (M\ttot N)^*\ot_k (M\ttot N)&\cong& (M^*\ot_k M)\sqot (N^*\ot_k N)\\
\nonumber (\psi\otimes_R\phi)\otimes_k(m\otimes_Rn)&\mapsto &(\phi\otimes_km)\sqot
(\psi\otimes_kn)
\end{eqnarray}

For any morphisms $f:X\lora Y, g:U\lora V$ in $\cat C$, by means
of \rref{rt6} we have the following diagram
\begin{equation}\label{rt7}\xymatrix{ 
(\fun F(Y)\ot_R \fun G(V))^*\otimes_k (\fun F(Y)\ot_R \fun G(V))
\ar[r]^\cong
&(\fun F(Y)^*\otimes_k \fun F(Y))\sqot(\fun G(Y)^*\otimes_k \fun G(V))
\ar[d]\\
(\fun F(Y)\ot_R \fun G(V))^*\otimes_k (\fun F(X)\ot_R \fun G(U))
\ar[u]_{\id\otimes_k\fun F(f)\otimes\fun G(g)}
\ar[d]_{(\fun F(f)\otimes\fun G(g))^*\otimes_k\id}&
L_0\sqot K_0\\
(\fun F(X)\ot_R \fun G(U))^*\otimes_k (\fun F(X)\ot_R \fun G(U))
\ar[r]_\cong
&(\fun F(Y)^*\otimes_k \fun F(Y))\sqot(\fun G(Y)^*\otimes_k \fun G(V))
\ar[u]
 }\end{equation}
Using the right exactness of the tensor product we see that 
$L\sqot K$ is the maximum quotient of $L_0\ot K_0$ that makes all the above
diagrams commutative. The claim of the proposition follows.\eee
\remark One can easily generalize the above proposition for more functors.\bigskip

Assume now that $\fun F:\cat C\to \RBimod$ is a monoidal functor, which means there
exists an $R$-bilinear natural isomorphism
$$\theta_{X,Y}:\fun F(X)\otimes_R\fun F(Y)\to \fun F(X\otimes Y)$$
satisfying the following identity (we assume for simplicity that $\cat C$ is strict, i.e. the
structure morphisms are identity morphisms)
\begin{equation}\label{rt9}\xymatrix{
\fun F(X)\otimes_R\fun F(Y)\otimes_R\fun F(Z)
\ar[rr]^{\theta_{X,Y}\otimes_R\id_{\fun F(Z)}}
\ar[d]_{\id_{\fun F(X)}\otimes_R\theta_{Y,Z}}
&&\fun F(X\otimes Y)\otimes_R\fun F(Z)\ar[d]_{\theta_{X\otimes Y,Z}}\\
\fun F(X)\otimes_R\fun F(Y\otimes Z)\ar[rr]_{\theta_{X,Y\otimes Z}}&&
\fun F(X\otimes Y\otimes Z)}
\end{equation} 
and there exists an isomorphism $\eta:\fun F(I)\to R$ ($I$ denotes the unit
object in $\cat C$) satisfying
\begin{equation}\label{rt10}\theta_{I,X}=\eta\otimes_R\id_X,\quad 
\theta_{X,I}=\id_X\otimes_R\eta
\end{equation}

It follows easily from definition that monoidal functors preserve rigidity. In fact we can
always choose $\ev_{\fun F(X)}$ and $\db_{\fun F(X)}$ to be $\fun F(\ev_X)$
and $\fun F(\db_X)$, respectively, in case $X$ is (left) rigid.

\begin{thm}\label{thm_rt4} Let $\cat C$ be a (strict) monoidal category and
$\fun F:\cat C\lora\RBimod$ be a monoidal functor with image in the subcategory of
left rigid bimodules. Let $L=\coend( \fun F)$. Then $L$ is a bialgebroid. If $\cat C$ is left 
rigid then $L$ is a Hopf algebroid.  If $\cat C$ is 
right rigid then $L$ is an opposite Hopf algebroid.\end{thm}
\proof We first show that $L$ is an $R$-bialgebroid. 
The product on $L$ is defined as follows. Consider the natural transformation
$$ \fun F(X)\ot_R \fun F(Y)\lora \fun F(X\ot Y)\lora \fun F(X\ot Y)\ot_R L\lora
\fun F(X)\ot_R \fun F(Y)\ot_R L$$
According to Proposition \ref{pro_rt3}, this natural transformation corresponds to a
morphism $m:L\sqot L\to L$, which according to Lemma \ref{lem_rt2} is a morphism
of $R$-coalgebroids. In other words, by means of the diagram in \rref{rt7}, $m$ is the 
unique map $L\sqot L\to L$ which satisfies the following diagram for all $X,Y\in\cat C$:
\begin{equation}\label{rt11}\xymatrix{
(\fun F(X)\otimes_R\fun F(Y))^*\otimes_k(\fun F(X)\otimes_R\fun F(Y))\ar[r]
\ar[d]_\cong&L\\
(\fun F(X)^*\otimes_k\fun F(X))\sqot(\fun F(Y)^*\otimes_R\fun F(Y))\ar[r]&
L\sqot L\ar[u]_m
}
\end{equation}
Further, since $\fun F(I)\cong R$, $R$ is a comodule over $L$. The coaction $R\to R\otimes_RL$
yields a morphism of $R$-coalgebroids $u:R\otimes_k R\to L$. It is easy to deduce
from Equations \rref{rt9}, \rref{rt10} and the universal property of $L$ the associativity
of $m$ and the unital property of $u$. Thus $L$ is an $R$-bialgebroid.

Assume that $\cat C$ is left rigid. We shall construct the antipode.
Recall that $L$ is a quotient of $L_0$, which is the direct sum of $\fun F(X)^*\otimes_R\fun F(X)$.
Set $M:=\fun F(X)$. For an element $\varphi\otimes_km$ of $M^*\otimes M$
we shall use the same notation to denote its image in $L$. Next, recall that the
defining relations for $L$ are obtained from morphism in $\cat C$.
In particular we deduce from the canonical morphism $\ev_X:X^*\otimes X\to I$
the following relation on $L$. Notice that 
$$\ev_X{}^*:I\to (X^*\otimes X)^{*}\cong X^*\otimes X^{**}$$
is nothing but $\db_{X^*}:I\to X^*\otimes X^{**}$.
By means of \rref{rt3} for the morphism
$\ev_X$ and using \rref{rt6}, we have the following commutative diagram, where $M:=\fun F(X)$,
\begin{equation}\label{rt12}\xymatrix{
R\otimes_k(M^*\otimes_kM)\ar[r]^{\id\otimes\ev_M}\ar[d]_{\ev_M{}*\otimes\id}&
R\otimes_kR\ar[d]^{(s\otimes t)}\\
(M^*\otimes_RM^{**})\otimes_k(M^*\otimes_R)\ar[r]\ar[d]_\cong&L\\
(M^{**}\otimes_kM^*)\sqot(M^*\otimes_kM)\ar[r]&L\sqot L\ar[u]_m
}
\end{equation}
Let $\{\varphi_j\}$, $\{\eta^j\}$ be dual bases with respect to $\db_{M^*}$, that is
$\db_{M^*}=\sum_j\varphi_j\otimes\eta^j\in M^*\otimes M^{**}$
(see Subsection \ref{sec_example}). Then \rref{rt12} amounts to the following relation
\begin{equation}\label{rt13}\sum_j (\eta^j\otimes_k\varphi)\circ(\varphi_j\otimes_km)=
\tau(\varphi(m))1
\end{equation}
where $1$ denotes the unit element in $L $ and $\circ$ denotes the product on $L$.
Similarly, by using the morphism $\db_X:I\to X\otimes X^*$ we obtain the following
relation on $L$:
\begin{equation}\label{rt14}(\varphi\otimes_k m_i)\circ(\eta\otimes_k\varphi^i)=
\si(\eta(\varphi))1 \end{equation}
where $\{m_i\}$ and $\{\varphi^i\}$ are dual bases with respect to the map
$\db_M:R\to M\otimes_RM^*$ ($M=\fun F(X))$), $\varphi\in M^*$, $\eta\in M^{**}$.

We define now the antipode $\nabla$. Recall that the map $\db_{k,M}:k\to M\otimes_RM^*$
was defined in Subsection \ref{sec_example} by $\db_{k,M}(1)=\sum_im_i\otimes\varphi^i$.
Define the map $\nabla_X$ for $M=\fun F(X)$
\begin{equation}\label{rt15}
\xymatrix{ M^*\ot_k M\ar_{\id\otimes\db_{k,M^*}\otimes\id}[d]   \ar[r]^{\nabla_X}  &L\ssot L\\
M^*\ot_kM^*\ot_RM^{**}\ot_kM\ar[r]_\cong
&(M^{**}\ot_k M^*)\ssot (M^*\ot_k M)\ar[u]
 }\end{equation}
where, in the tensor product $M^{**}\ot_kM^*$, we use the convention
 that the action on $M^*=\fun F(X^*)$ is denote by
$\tau$ and the action on $M^{**}$ is denoted by $\si$. 
It is straightforward to check the commutativity of the following diagram
$$\xymatrix{\fun F(Y)^*\otimes_k\fun F(X)\ar[r]^{\fun F(f)^*\otimes\id}
\ar[d]^{\id\otimes\fun F(f)}&
\fun F(X)^*\otimes_k\fun F(X)\ar[d]^{\nabla_X}\\
\fun F(Y)^*\otimes_k\fun F(Y)\ar[r]_{\nabla_Y}&L\ssot L
}$$
Thus the universal property of $L$ yields a morphism $\nabla:L\to L\ssot L$ which
we will show to be the antipode of $L$.
Explicitly we have
\begin{equation}\label{rt16}\nabla(\varphi\otimes m)=
(\eta^j\otimes_k \varphi)\ssot (\varphi_j\otimes_km)
\end{equation}
where the dual bases $\{\eta^j,\varphi_j\}$ are defined
above.
Now, the equation \rref{eq341}, \rref{eq342} for $\nabla$ can be easily deduced from
\rref{rt13}, \rref{rt14}. Let us show \rref{eq341} for $h=\varphi\otimes_k m$. The
left hand side of \rref{eq341} is equal to
$$(\eta^j\otimes_k \varphi)\circ(\varphi_j\otimes_k m_i\tsot (\varphi^i\otimes_k m))
=\varphi(m_i)(\varphi^i\otimes_k m)=\varphi\otimes_km$$
where in the first equation we used \rref{rt14}. We thus showed that $L$ is an
$R$-Hopf algebroid.

If  $\cat C$ is right rigid (in this case the image of $\fun F$ lies in the subcategory
of rigid bimodules), the opposite antipode is induced from the maps
$$\xymatrix{ M^*\ot_k M\ar_{\id\otimes\db_{^*M}\otimes\id}[d]\ar^{\bar\nabla_M}[r]&L\ttot L
\\
 M^*\ot_k{}^*M\ot_RM\ot_kM\ar[r]_\cong&
  (M\ot_k {}^*M)\ttot (M^*\ot_k M)\ar[u]
  }$$
where $M=\fun F(X)$, $X\in\cat C$.
\eee

\begin{cor}\label{repr-bialg} Let $\cat C$ be a small locally finite
$k$-linear abelian monoidal category and $\fun F:\cat C\lora
\RBimod$ be a faithful exact, monoidal functor with image in the subcategory
of right f.g. projective modules. Let $L=\coend (\fun F)$. Then
$L$ is semi-transitive coring with respect to the actions $(\tau,\si)$ and
$\fun F$ induces a monoidal equivalence between $\cat C$
and $\comodL$. Conversely, let $L$ be a  bialgebroid, semi-transitive as a coring
with respect to the actions $(\tau,\si)$, and
$\fun F$ be the forgetful functor into the category of $R$-bimodules.
Then $\fun F$ has image in the subcategory of left rigid bimodules and
$L\cong\coend(\fun F)$.
\end{cor}
\proof  We first notice that if $\fun F:\cat C\to \cat D$ is at the same time
a monoidal functor and an equivalence then $\fun F$ is a monoidal
equivalence (i.e. the quasi-inverse to $\fun F$ is also a monoidal functor).
Indeed, let $\theta$ and $\eta$ be the structure morphism for $\fun F$ as in
\rref{rt9}, \rref{rt10}. By definition of the quasi-inverse we have the natural isormophims
$\fun{F}\fun G(U)\cong U$ and $\fun{GF}(X)\cong X$, cf. \cite[Section IV.4]{maclane}.
Then (still assuming the categories to be strict for
simplicity) we define the monoidal functor structure for the quasi-inverse
$\fun G$ of $\fun F$ as follows:
\begin{eqnarray}\label{rt17}&
\zeta_{U,V}:\fun G(U)\otimes \fun G(V)\cong
\fun{G}\fun F(\fun G(U)\otimes\fun G(V))\stackrel{\fun G(\theta^{-1})}
\lora \fun G(U\otimes V)\\
& \xi:\fun G(I_{\cat D})\stackrel{\eta^{-1}}\lora\fun {GF}(I_{\cat C})
\cong I_{\cat C}\end{eqnarray}

Assume that we have $\fun F:\cat C\lora \RBimod$ as required.
Let $L$ be the $\coend$   of $\fun F$. Then, by virtue of Theorem \ref{thm_deligne},
$L$ is semi-transitive as an $R$-coring with respect to the pair of actions  $\si,\tau$
and the induced functor $\bar{\fun F}$ is an equivalence of abelian categories.

On the other hand, by virtue of Theorem \ref{thm_rt4}, $L$ is an $R$-bialgebroid
and $\bar{\fun F}$ is a monoidal functor to the category of $L$-comodules,
 thus $\bar{\fun F}$ is
 a monoidal equivalence.

Assume that the bialgebroid $L$ is  semi-transitive coring with respect to the
actions $(\tau,\si)$. Then the forgetful functor has image
in the subcategory of rigid bimodules. Let $L'$ be the $\coend$   of this functor, we have
a morphism of bialgebroids $L'\lora L$ which is an isomorphism of corings, by
virtue of Theorem \ref{thm_deligne}, hence $L'\cong L$ as bialgebroids.\eee

In what follows we will consider only Hopf algebroids with opposite antipode.
\definition ({\em Semi-transitive Hopf algebroids}) A Hopf algebroid $H$
is said to be semi-transitive over $R$ if the following conditions are satisfied:
\begin{enumerate}\item $H$ is semi-transitive  as an $R$-coring
with respect to the actions $(\si,\tau)$.
\item an $H$-comodule is left rigid as an $R$-bimodule iff it is right rigid.\end{enumerate}

\begin{thm}\label{thm_repr1}Let $\cat C$ be a small locally finite
$k$-linear abelian rigid monoidal category and $\fun F:\cat C\lora
\RBimod$ be a faithful exact, monoidal functor. Let $H=\coend (\fun F)$.
Then $H$ is a semi-transitive Hopf algebroid (with opposite antipode) and $\fun F$ induces a
monoidal equivalence between $\cat C$ and finitely generated (over $R$)
right $H$-comodules. Conversely, let $H$ be a semi-transitive Hopf algebroid and
$\fun F$ be the forgetful functor from the category of finitely generated (over $R$)
$H$-comodules to $\RBimod$. Then $H\cong \coend \fun F$.
\end{thm}
\proof $H=\coend(\fun F)$ is obviously a Hopf algebroid with opposite antipode.
The equivalence is established by the corollary above. Also, from the construction, we see
that $H$ is a semi-transitive as a coring with respect to the actions $(\si,\tau)$.

 It remains to show that an $H$-comodule is left rigid
if and only if it is right rigid. Let $M$ be a right $H$-comodule
which left rigid as $R$-bimodule. Then by the equivalence, $M\cong
\fun F(X)$ for a certain $X\in \cat C$. Hence $M$ is rigid for $X$
is rigid. Conversely, if $M$ is  a right  rigid $R$-bimodule, then
${}^*M$ is left rigid. Thanks the opposite antipode ${}^*M$ has a
structure of $H$-comodule, hence ${}^*M\cong \fun F(Y)$ for a
certain $Y$ in $\cat C$. Therefore $M\cong \fun F(Y^*)$; hence
rigid.

Now assume that $H$ is a semi-transitive Hopf algebroid. Then the forgetful functor
$\fun F:\comod$-$H\lora \RBimod$ has image
in the subcategory of rigid bimodules. Since $H$ is semi-transitive, this functor
is exact (and obviously faithful being forgetful functor). Thus, we can reconstruct the
$\coend$ of this functor. By virtue of Theorem \ref{thm_deligne},
$\coend(F)\cong H$ as corings;
hence they are isomorphic as Hopf algebroids.\eee

\remark The condition (ii) in the definition of semi-transitive Hopf algebroid is not
natural. In fact, it is used only for the formulation of Theorem \ref{thm_repr1}. In other
words, Theorem \ref{thm_repr1} states that one can ``fully'' reconstruct a Hopf algebroid
from a faithful, exact monoidal functor, in the sense that if we repeat this process we
will obtain the same Hopf algebroid. However, we do not have a good criterion for a
Hopf algebroid to be reconstructible from its category of comodules. The reader is also
referred to \cite{brug1} for some problems related to the notion of transitivity.\\

Theorem \ref{thm_repr1} has  an interesting consequence on characterizing abstract
rigid monoidal categories. First, we mention a result of \cite{ph00}

{\em Let $\cat C$ be a small abelian rigid monoidal
 category. Then there exists an exact faithful monoidal functor $\cat C\lora \RBimod$
for a certain ring $R$.}

By using the above result of reconstruction and representation, we can easily
deduce the following result
\begin{cor}\label{cor_repr} Let $\cat C$ be a small $k$-linear
locally finite abelian rigid monoidal category. Then there exists a ring $R$ such 
that $\cat C$ is monoidally equivalent with the category of f.g. projective 
$R$-comodules over a certain semi-transitive Hopf algebroid over $R$.\end{cor}
\bigskip

\centerline{\bf Acknowledgment}
\bigskip

 This work is supported by the National Program for Basic Sciences Research, Vietnam.
 A part of this work was carried out during the author's visit at the ICTP, Trieste, Italy,
to which he would like to express his sincere thank for providing excellent working condition and
financial support.
 The author also thanks Professors Nguyen Dinh Cong and Do Ngoc Diep for
stimulating discussions.
Finally he would like to thank the referee for carefully reading the manuscript, pointing out misprints and
making helpful remarks, comments which substantially improved the manuscript.

\end{document}